# The Algebra and Geometry of Isometries

**Dedicated to Catherine Schmidt**

**John E. Connett, PhD**
**Division of Biostatistics**
**School of Public Health**
**University of Minnesota**

**john-c@ccbr.umn.edu**

**May 23, 2014**

**1. On the Pitcher's Mound.** The pitcher holds the ball behind his back in his glove. He waits for the catcher's signal. He rears back and throws a slow-spinning curve ball. The batter hits a low liner to shallow center field. The center fielder efficiently fields the ball on the second hop, and throws it to first base, where it is caught for the first out, The first baseman relays the ball back to the pitcher. The pitcher returns the ball to the same position in his glove, behind his back, as he gets ready for the next pitch (Figure 1).

Unnoticed by pitcher, fielders, the next batter and fans in the stadium[1], there are two points on the surface of the ball which ends up exactly where it started before all the motions of the ball occurred. The surface of the ball is approximately a 2-sphere; thus, unless the batter hit the ball with such force that the cover was torn off, the sequence of motions has two fixed points.

Through the various motions of throwing, spinning, being hit by the bat, bouncing on the field, fielded and thrown to first base and then back to the pitcher, the points on the ball have undergone an **isometry** - that is, the points have been subjected to a sequence of motions that preserve distances between all pairs of points. Because the isometry of the ball resulted from a sequence of physical motions, it is orientation-preserving (defined below). The objective of this paper is to examine isometries on the 2-sphere and in the plane, and to compare algebraic and geometric approaches for determining the fixed axes (for the 2-sphere) and fixed points (for the plane) of orientation-preserving isometries. A good general reference on isometries in the plane is [Pedoe (1970)].

**2. Isometries in the Plane, $R^2$:** There are three types of isometries in the plane: translations, rotations around a fixed point, and reflections through a line.

A *translation* $T_V : R^2 \to R^2$ moves each point $X$ in the plane to $X + V$, where $V$ is a fixed vector.

A *rotation* $A_{P,\theta} : R^2 \to R^2$ around a point $P$ and through a counterclockwise angle $\theta$ is defined by

$$A_{P,\theta}(X) = P + R_\theta(X - P),$$

where $R_\theta$ is a **central** rotation, representing counterclockwise rotation through angle $\theta$ around the origin $O$. Note that $R_\theta$ is a linear function which can be represented by a $2 \times 2$ matrix $M_\theta$,

$$M_\theta = \begin{bmatrix} \cos(\theta) & -\sin(\theta) \\ \sin(\theta) & \cos(\theta) \end{bmatrix}.$$

The point $P$ is the fixed point of $A_{P,\theta}$, referred to in the following as a *pivot* point.

---
[1] And possibly of no interest to anyone present …

A *reflection* $L_\ell$ through the line $\ell$ transforms the point $X$ into the mirror image $X' = L_\ell(X)$. A reflection through the line $\ell$ can be written as the composition of a reflection through the $X$ axis followed by a rotation $R_\theta$ around the origin, where $\theta$ is the angle which $\ell$ makes with the $X$ axis.

Facts:

1. The set of all isometries of the plane forms a group.

2. The composition of two reflections (across lines $\ell$ and $\ell'$) in general is a rotation. If the angle between the lines is $\alpha$, then the angle of the rotation is $2\alpha$: see Figure 2. If the lines are parallel, the composition of the two reflections is a translation in a direction perpendicular to the lines, and of magnitude equal to twice the distance between the lines.

3. The composition of two rotations in general is another rotation. If the corresponding angles $\theta$ and $\theta'$ sum to zero, the composition is a translation.

4. The group of isometries is nonabelian. In general, two reflections do not commute; a reflection and a rotation do not commute; two rotations do not commute; a translation and a reflection do not commute; a translation and a rotation do not commute. Two translations do commute, and two rotations pivoting around the same point commute. Two reflections commute if their lines are parallel.

5. Orientation: An isometry $S$ which preserves orientation has the property that, if a triangle has vertices in the order $A, B, C$ when the triangle is traversed in a counterclockwise direction, then for the image of that triangle under the transformation $S$, the vertices $A', B', C' = S(A), S(B), S(C)$ are also in counterclockwise order. Rotations and translations preserve orientation; reflections do not.

The main objective in the following is to compare and contrast algebraic and geometric/constructive approaches to some aspects of isometries of the plane. We will present this as a series of problems:

**Problem P1:** Given a line segment $XY$ in the plane and another segment $X'Y'$ which you know to be its image under an orientation-preserving isometry, find the isometry.

First, note that $length(XY) = length(X'Y')$. Next, if $XY$ and $X'Y'$ are parallel and the direction from $X'$ to $Y'$ is the same as the direction from $X$ to $Y$, then the unknown isometry is a translation $T_V$, where $V = X' - X = Y' - Y$.

If $XY$ and $X'Y'$ are not parallel, then the isometry is a rotation $A_{P,\theta}$ with pivot $P$ and angle $\theta$, and the task is to find $P$ and $\theta$.

**Algebraic approach:**

In the following, points in the plane are thought of as $2 \times 1$ column vectors, but for convenience these can be represented as the transpose of $1 \times 2$ row vectors, i.e.,

$$\begin{pmatrix} a \\ b \end{pmatrix} = (a,b)^T.$$

Let $X = (x_1, x_2)^T$, $Y = (y_1, y_2)^T$, and $X' = (x'_1, x'_2)^T$, $Y' = (y'_1, y'_2)^T$, and assume $A_{P,\theta}$ is a rotation which takes $XY$ to $X'Y'$. From the definition of rotation, this means that

$P + R_\theta(X - P) = X'$ and $P + R_\theta(Y - P) = Y'$. Since $R_\theta$ is a linear function, this implies that

$P - R_\theta(P) + R_\theta(X) = X'$ and $P - R_\theta(P) + R_\theta(Y) = Y'$, and further that

$R_\theta(X - Y) = X' - Y'$.

Therefore

$\begin{bmatrix} \cos(\theta) & -\sin(\theta) \\ \sin(\theta) & \cos(\theta) \end{bmatrix} \cdot \begin{bmatrix} x_1 - y_1 \\ x_2 - y_2 \end{bmatrix} = \begin{bmatrix} x'_1 - y'_1 \\ x'_2 - y'_2 \end{bmatrix}$, which is equivalent to the equations

$\cos(\theta) \cdot (x_1 - y_1) - \sin(\theta) \cdot (x_2 - y_2) = (x'_1 - y'_1)$ and

$\sin(\theta) \cdot (x_1 - y_1) + \cos(\theta) \cdot (x_2 - y_2) = (x'_2 - y'_2)$, which can be reconfigured as

$\begin{bmatrix} x_1 - y_1 & -(x_2 - y_2) \\ (x_2 - y_2) & (x_1 - y_1) \end{bmatrix} \begin{bmatrix} \cos(\theta) \\ \sin(\theta) \end{bmatrix} = \begin{bmatrix} x'_1 - y'_1 \\ x'_2 - y'_2 \end{bmatrix}$, yielding

$\begin{bmatrix} \cos(\theta) \\ \sin(\theta) \end{bmatrix} = \begin{bmatrix} x_1 - y_1 & -(x_2 - y_2) \\ (x_2 - y_2) & (x_1 - y_1) \end{bmatrix}^{-1} \cdot \begin{bmatrix} x'_1 - y'_1 \\ x'_2 - y'_2 \end{bmatrix}.$

It is straightforward to check, using the fact that $length(XY) = length(X'Y')$, that this linear-algebraic solution has the property that $\cos^2(\theta) + \sin^2(\theta) = 1$.

Thus $\theta = \arccos(\cos(\theta))$, where the interior "$\cos(\theta)$" is regarded as a solution to the linear equations above. The pivot point $P$ can be derived as follows:

$$P - R_\theta(P) = X' - R_\theta(X), \text{ so}$$

$$P = (I - R_\theta)^{-1} \cdot (X' - R_\theta(X)),$$

where $I$ is the identity transformation.

**Geometric approach:**

See Figure 3. The unknown isometry rotates by the angle $\theta$ around a point $P$, taking $X$ to $X'$ and $Y$ to $Y'$. Therefore P must lie on the perpendicular bisector $\ell_X$ of $XX'$, and similarly on the perpendicular bisector $\ell_Y$ of $YY'$. Therefore $P$ can be constructed as the point of intersection of these two lines.

Note that in Figure 3, the triangles $\Delta PXY$ and $\Delta PX'Y'$ are congruent, since the three corresponding sides are equal. From this it follows that $\angle XPX' = \angle YPY'$. This common angle is therefore the angle $\theta$ of the rotation, and the rotation $A_{P,\theta}$ is completely defined.

**Problem P2:** Given two rotations $A_{G,\alpha}$ and $A_{H,\beta}$, find the pivot point and the angle of rotation of the composition $A_{G,\alpha} \circ A_{H,\beta}$.

We will assume in the following that $\beta \neq -\alpha$, since otherwise the composition is trivially shown to be a translation $T_V$, where $V = G + H$.

**The "Bird's Eye" Proof:**

As a preliminary to both the algebraic or geometric approaches we will give an intuitive proof that the angle of rotation of the composite isometry must be $\theta = \alpha + \beta$.

Assume a bird is flying very high above the plane, so high that the points $G, H$ and the unknown pivot point $P$ of the composite cannot be distinguished – they all look like one point. The rotation $A_{H,\beta}$ rotates the entire plane through the angle $\beta$. The rotation $A_{G,\alpha}$ then rotates the entire plane through the angle $\alpha$. To the bird, the result of the composite appears to be a rotation of the whole plane through the angle $\beta + \alpha$, around a single point, which because of its great altitude, the bird cannot distinguish the pivots as separate points.. Therefore the angle of rotation of the composite isometry must be $\alpha + \beta$.

**Algebraic Approach to Finding P:**

As above, the two isometries are $A_{H,\beta}$ and $A_{G,\alpha}$, and the composition is $A_{P,\theta} = A_{G,\alpha} \circ A_{H,\beta}$. Since $A_{H,\beta}(X) = H + R_\beta(X - H)$, we have that

$$A_{G,\alpha} \circ A_{H,\beta}(X) = G + R_\alpha(H + R_\beta(X - H) - G).$$

Because $R_\alpha$ and $R_\beta$ are linear functions, this can be written as

$$G + R_\alpha(H) + R_\alpha(R_\beta(X)) - R_\alpha(R_\beta(H)) - R_\alpha(G)$$
$$= G + R_\alpha(H) + R_{\alpha+\beta}(X) - R_{\alpha+\beta}(H) - R_\alpha(G).$$

The unknown isometry acting on $X$ is

$$A_{P,\theta}(X) = A_{P,\alpha+\beta}(X) = P + R_{\alpha+\beta}(X - P) = P + R_{\alpha+\beta}(X) - R_{\alpha+\beta}(P).$$

Therefore

$$P - R_{\alpha+\beta}(P) = G + R_\alpha(H) - R_{\alpha+\beta}(H) - R_\alpha(G), \text{ and}$$

$$(I - R_{\alpha+\beta}) \cdot P = G + R_\alpha(H) - R_{\alpha+\beta}(H) - R_\alpha(G), \text{ which yields the solution for } P:$$

$$P = (I - R_{\alpha+\beta})^{-1} \cdot (G + R_\alpha(H) - R_{\alpha+\beta}(H) - R_\alpha(G)).$$

A somewhat ugly expression to be sure, but easily computable in a matrix-language program like MATLAB or Mathematica or SAS PROC IML.

**Geometric Approach to Finding P:**

This is basically a variant of the geometric argument in Problem P1. Figure 4 shows a line segment $XY$, its image under the first isometry $A_{H,\beta}(XY) = X'Y'$, and the subsequent image under the composite isometry, $A_{P,\alpha+\beta}(XY) = A_{G,\alpha}(X'Y') = X''Y''$. In Figure 4, the first rotation is around the point $H = (1,0)^T$, through an angle of $\beta = \pi/2$ radians $= 90°$. The second rotation is around the point $G = (0,0)^T$ and the associated angle of rotation is $\alpha = \pi/4 = 45°$. As in Problem P1, the point $P$ is the intersection of the perpendicular bisectors $\ell_X$ and $\ell_Y$ of the line segments $XY$ and $X''Y''$. The angle of rotation from the "bird's eye" argument is $\alpha + \beta = 3\pi/4 = 135°$.

In this example, both the algebraic and geometric approaches give $P = (0.7071, 0.2929)^T$.

**Problem P3:** Given a rotation $A_{P,\theta}$, find two line-reflections whose composition equals $A_{P,\theta}$.

This is an easy problem, especially in view of Figure 2. There is not a unique solution. Choose any two lines through the point $P$ such that the angle between them is $\theta/2$.

## 3. Isometries of the 2-Sphere.

The surface of a baseball is approximately a 2-dimensional sphere, $S^2$. There are essentially 4 categories of isometries on $S^2$:

1. Rotations around a fixed axis.

2. Reflections through a great circle (i.e., an 'equator').

3. Reflections through the center of a ball which has the sphere as its surface.

4. Compositions of reflections and rotation.

In the following we will consider $S^2$ to be the 2-sphere of radius 1 centered around the origin in $R^3$.

There is no analogue on $S^2$ of the translations which occur in $R^2$. Only the isometries in the first category preserve orientation. As with $R^2$, the isometries of $S^2$ form a group. In fact the orientation-preserving isometries are identical to the subgroup $SO(3)$ of the linear isometries of 3-dimensional space.

In the plane, geodesic segments are straight lines, i.e. curves of least distance joining two points.. On the 2-sphere, the analogues are segments of great circles, i.e., circles which are the intersection of a plane through the origin $O = (0,0,0)^T$ in 3-space, and the 2-sphere. Given two points on the sphere, the shortest curve on the sphere which joins them is a segment of a great circle. Ordinarily this circle is unique; however if the points happen to be antipodal, there are infinitely many great circles on which both points lie.

**Problem S1:** Given two great-circle segments $XX'$ and $YY'$, where $X$ and $X'$ are non-antipodal points on the 2-sphere and the segments are of equal length, find the orientation-preserving isometry of $S^2$ which transforms $XX'$ to $YY'$.

**Geometric Approach:**

Refer to Figure 5. This is directly analogous to the geometric approach to the same problem in the plane. First, draw the great-circle segments $XX'$ and $YY'$. Let $\ell_X$ and $\ell_Y$ be the great circles which are the perpendicular bisectors of $XX'$ and $YY'$ respectively. Any two distinct great circles on the 2-sphere intersect in two points $P$ and $P'$ (which are antipodal to each other) and determine an axis through the origin. As in the argument for Problem P1 in the plane, $\Delta PXY$ and $\Delta PX'Y'$ are congruent spherical triangles. Thus rotation around the axis $PP'$ transforms $XY$ to $X'Y'$. The amount of this rotation is determined by measuring either of the angles $\angle XPX'$ or $\angle YPY'$ (see Figure 5). Note that both $P$ and $P'$ are fixed points of this rotation.

**Algebraic Approach:**

As with the geometric approach, the first step to finding the isometry that transforms $XY$ to $X'Y'$ is to find the axis $PP'$. First, note that the vector $X - X'$ is orthogonal to the disc bounded by the great circle which is the perpendicular bisector of the segment $XX'$. Similarly, $Y - Y'$ is orthogonal to the disc of the great circle that is the perpendicular bisector of $YY'$. These two discs intersect in a line through the origin. This line is orthogonal to both $X - X'$ and $Y - Y'$, which means that the direction of this line is determined by the cross product[2] of $X - X'$ and $Y - Y'$, i.e.,

$$U = (X - X') \times (Y - Y').$$

Normalizing this vector by dividing through by its length $|U|$ gives $P = U/|U|$, and the antipodal point is $P' = -P$.

Determining the angle of rotation can again be accomplished by computation of another cross product, namely, $V = X \times X'$, which has magnitude $|X| \cdot |X'| \cdot \sin(\theta)$, i.e.,
$\theta = \arcsin(|V|/(|X| \cdot |X'|))$.

**Problem S2:** Given two rotations of $S^2$, find the axis corresponding to their composition and the angle of rotation.

**Algebraic Approach:**

Let the two isometries be $A_{G,\alpha}$ and $A_{H,\beta}$, where $G$ and $H$ are points on the 2-sphere which determine the axes of rotation (i.e. one axis is determined by $G$ and its antipodal point $-G$, and similarly the other axis is determined by $H$ and $-H$. Both $A_{G,\alpha}$ and $A_{H,\beta}$ can be represented by matrices, that is,

---

[2] See the Appendix for a discussion of cross products, eigenvalues and eigenvectors.

$$A_{G,\alpha} \leftrightarrow M_{G,\alpha} = \begin{bmatrix} a_1 & a_2 & a_3 \\ b_1 & b_2 & b_3 \\ c_1 & c_2 & c_3 \end{bmatrix} \text{ and } A_{H,\beta} \leftrightarrow M_{H,\beta} = \begin{bmatrix} d_1 & d_2 & d_3 \\ e_1 & e_2 & e_3 \\ f_1 & f_2 & f_3 \end{bmatrix}.$$

The composition of $A_{H,\beta}$ followed by $A_{G,\alpha}$ can be represented by $M_{P,\theta} = M_{G,\alpha} \cdot M_{H,\beta}$ where $P$ and its antipodal point $P' = -P$ define the to-be-determined axis of rotation. All $3 \times 3$ matrices have at least one real eigenvalue (because the characteristic polynomial is a cubic; see [Axler, 1997]). Because $A_{P,\theta}$ is an isometry, the real eigenvalues of $M_{P,\theta}$ must be $+1$ or $-1$. In fact, for orientation-preserving isometries, the real eigenvalue must be $+1$. Here's why. The eigenvalues are roots of the characteristic polynomial. The determinant of $M_{P,\theta}$ is the product of the eigenvalues. Two of the eigenvalues are complex; designate these by $a + ib$ and $a - ib$. The product of these is $a^2 + b^2$. The determinant for an orientation-preserving transformation must be positive. Therefore the real eigenvalue must be $+1$. Thus, let $P$ be the eigenvector corresponding to the real eigenvalue[3].

Unlike composition of rotations in the plane, the angle of rotation $\theta$ of the composition $M_{P,\theta}$ is not equal to the sum of the angles of the two original isometries. Clearly the "bird's eye" proof for the plane does not apply to the 2-sphere. Here is an example:

Let $A_{H,\beta}$ be defined as a rotation around the $z$ axis in $R^3$ of $\pi/6 = 0.5236$ radians, and $A_{G,\alpha}$ is rotation around the $y$ axis of $\pi/4 = 0.7854$ radians. The pivot points of $A_{H,\beta}$ and $A_{G,\alpha}$ respectively are $H = (0,0,1)^T$ and $G = (0,1,0)^T$. The pivot point $P$ of $A_{P,\theta}$ is the eigenvector $P = (-0.2195, -0.8192, 0.5299)^T$ and the angle of rotation around the axis $PP'$, as given by the cross product approach in Problem S1, is $\theta = 0.9363 \neq 1.3090 = \alpha + \beta$. The angle of rotation can also be found by consideration of the two complex eigenvalues of $M_{P,\theta}$. In this case, those two eigenvalues are:

$$e_1 = .5927 + i*.8054 \text{ and } e_2 = .5927 - i*.8054,$$

which implies that the angle of rotation around the axis $PP'$ is $\theta = \arccos(0.5927) = 0.9363$ radians.

**Geometric Approach:**

The geometric approach in this case is exactly analogous to Problem P2 in the plane: Fix a segment $XY$ on the sphere, apply the first isometry $A_{H,\beta}$ to obtain the image $X'Y'$ and then

---

[3] See the Appendix for discussion of cross products, eigenvalues and eigenvectors

apply $A_{G,\alpha}$ to yield $X''Y''$. Then use the geometric approach above for Problem S1 for the 2-sphere to construct $P$ and $\theta$.

## 4. How to Find Fixed Points on the Baseball:

If the pitcher is overcome by curiosity to find what the axis of rotation and fixed points are for the baseball after it has gone through its various spins, throws, and bounces, he can do the following. Before the first pitch, mark two non-antipodal points $X$ and $Y$ and connect them by the segment of the great circle on the sphere determined by these two points. The pitcher can then photograph the sphere, being sure to include the segment $XY$ in the picture, using his cell-phone. Then he throws the ball and eventually it is returned to him. He photographs the ball again, noting the position of the transformed segment $X'Y'$. Then he applies either the geometric or algebraic approach to Problem S1 for the 2-sphere. Hopefully the opposing side will not call for a delay-of-game penalty!

**Appendix: Cross-Products, Eigenvectors, and Eigenvalues (see Axler[2]).**

In $R^3$, the **cross-product** $X \times Y$ of two vectors $X = (x_1, x_2, x_3)^T$ and $Y = (y_1, y_2, y_3)^T$ is defined as

$$X \times Y = (x_2 y_3 - x_3 y_2, -x_1 y_3 + x_3 y_1, x_1 y_2 - x_2 y_1)^T.$$

Basic properties of the cross product are (1) $X \times Y$ is orthogonal to both $X$ and $Y$, and (2) the length $|X \times Y|$ of $X \times Y$ is equal to $|X| \cdot |Y| \cdot \sin(\theta)$, where $\theta$ is the angle between $X$ and $Y$.

If $A: R^n \to R^n$ is a linear transformation, then an **eigenvector** $X$ of $A$ is a vector such that

$A(X) = \lambda \cdot X$, where $\lambda$ is a real number, called an **eigenvalue**.

An essential property of eigenvalues is that they are roots of the **characteristic polynomial** of $A$, where the characteristic polynomial in $\lambda$ is defined as the determinant of $\lambda \cdot I - A$, i.e.,

$$\det(\lambda \cdot I - A),$$

where $I$ is the $n \times n$ identity matrix.

This polynomial in $\lambda$ has degree $n$. In general, the roots of this polynomial may be real or complex; only the real roots correspond to actual eigenvectors of $A$. Eigenvectors are not unique; any vector with the same *direction* as a given eigenvector is also an eigenvector with the same eigenvalue. The essential fact about eigenvectors of a linear transformation is that the transformation *preserves the original direction* of the eigenvector.

All graphics (except Figure 1) and computations in this paper were done using SAS/GRAPH and SAS PROC IML [SAS9.2, 2011].

**Figure 1:**

**The Baseball Just Before the First and Second Pitches**

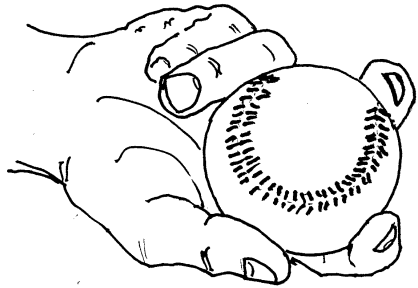 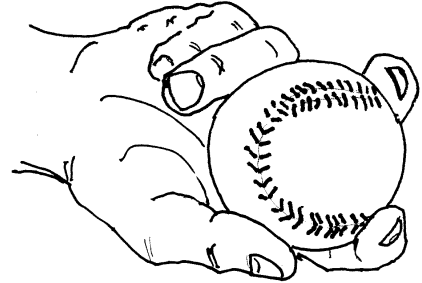

**First Pitch**                      **Second Pitch**

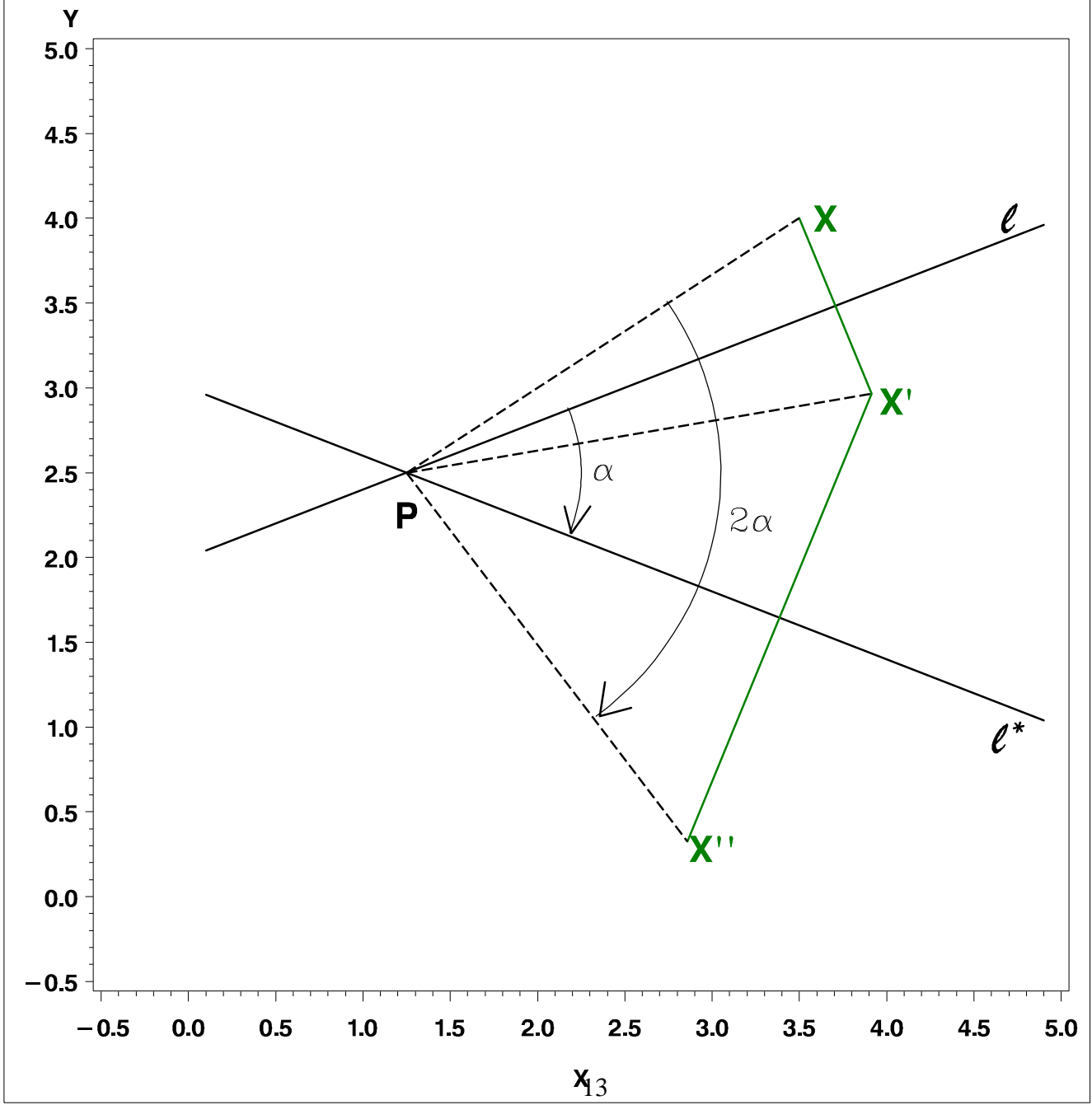

Figure 2: Two Line–Reflections Equal a Rotation

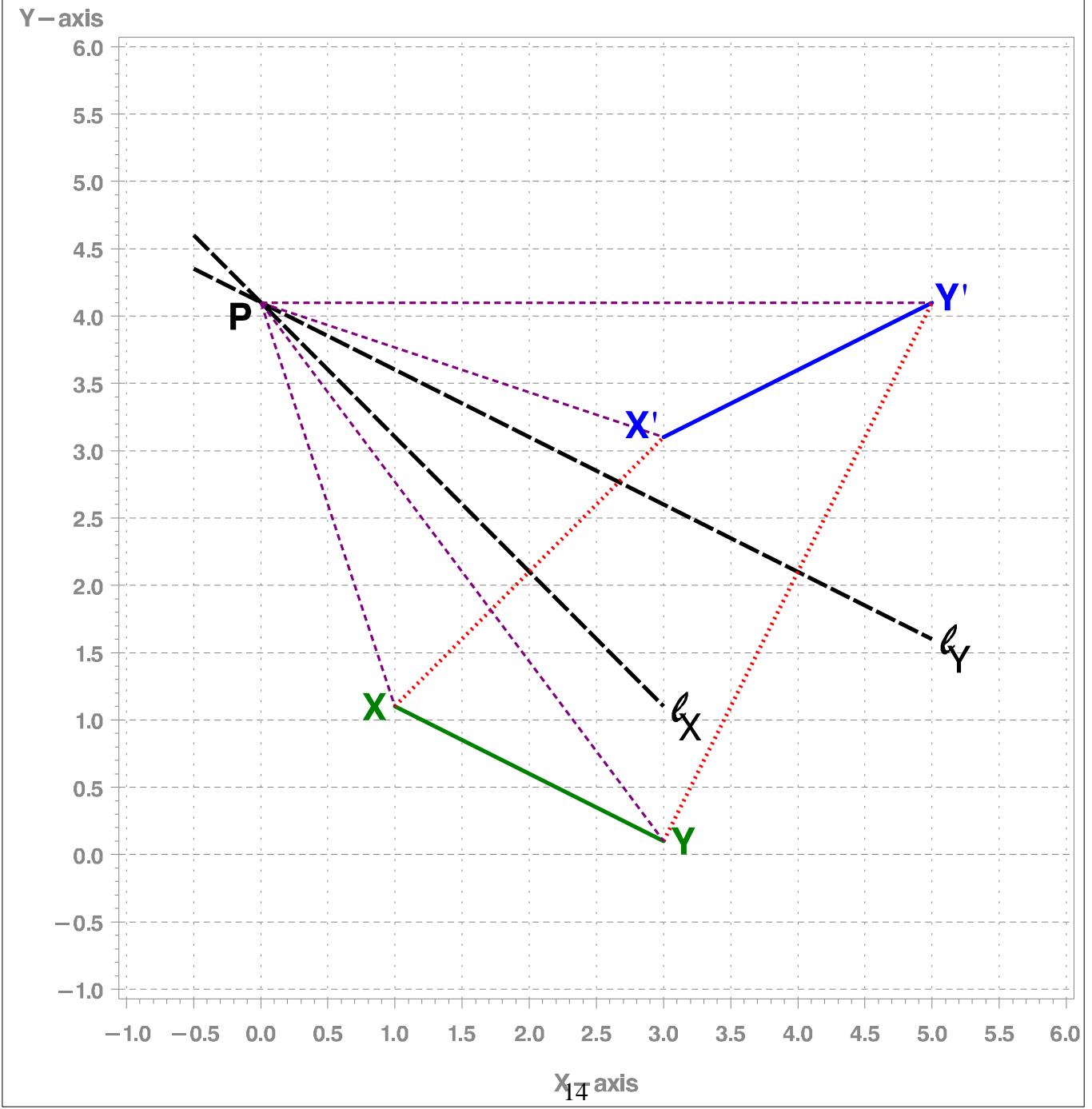

Figure 3: Locating the Pivot Point P, and Finding the Angle of Rotation, Given Line Segments XY and X'Y'

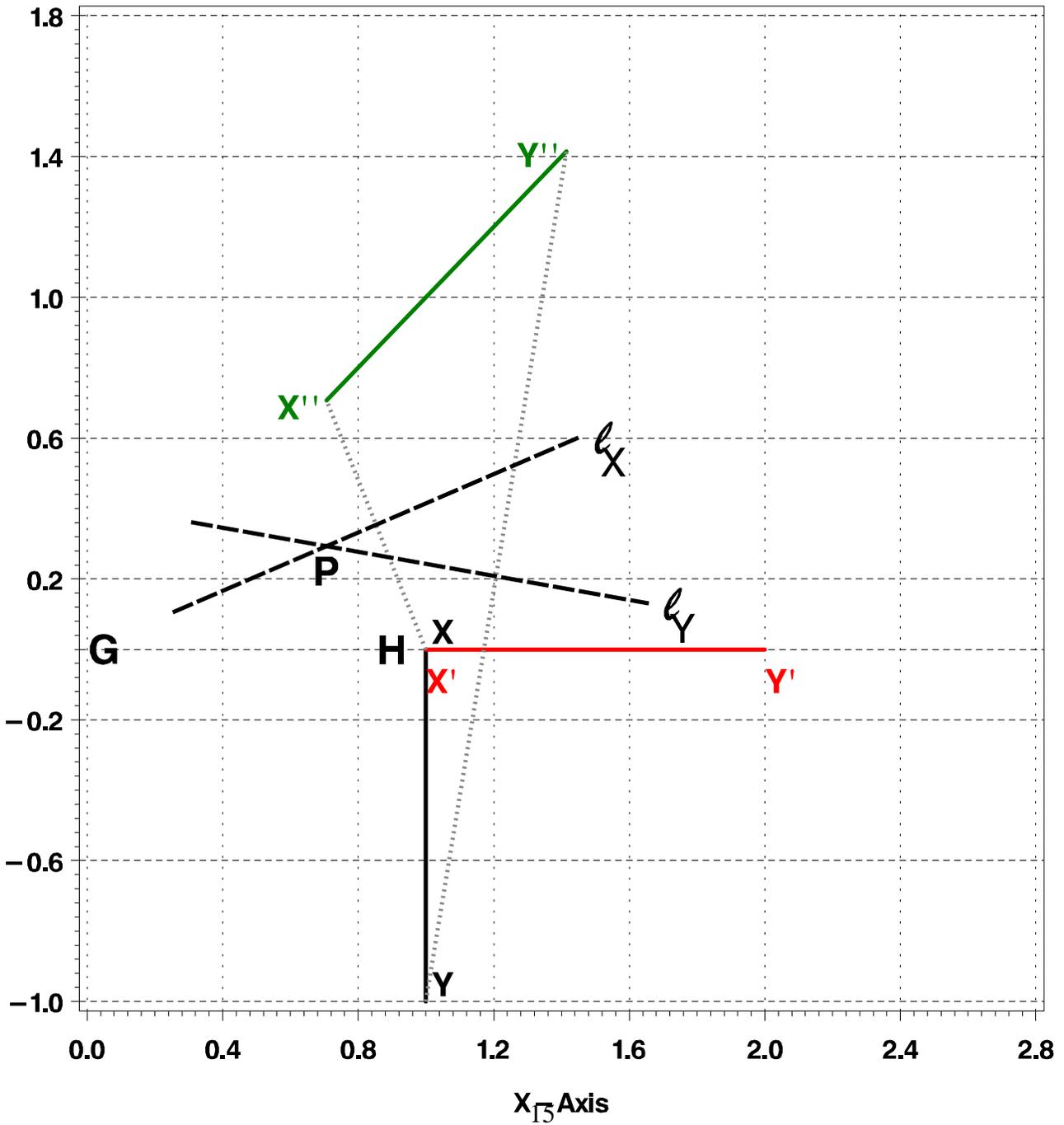

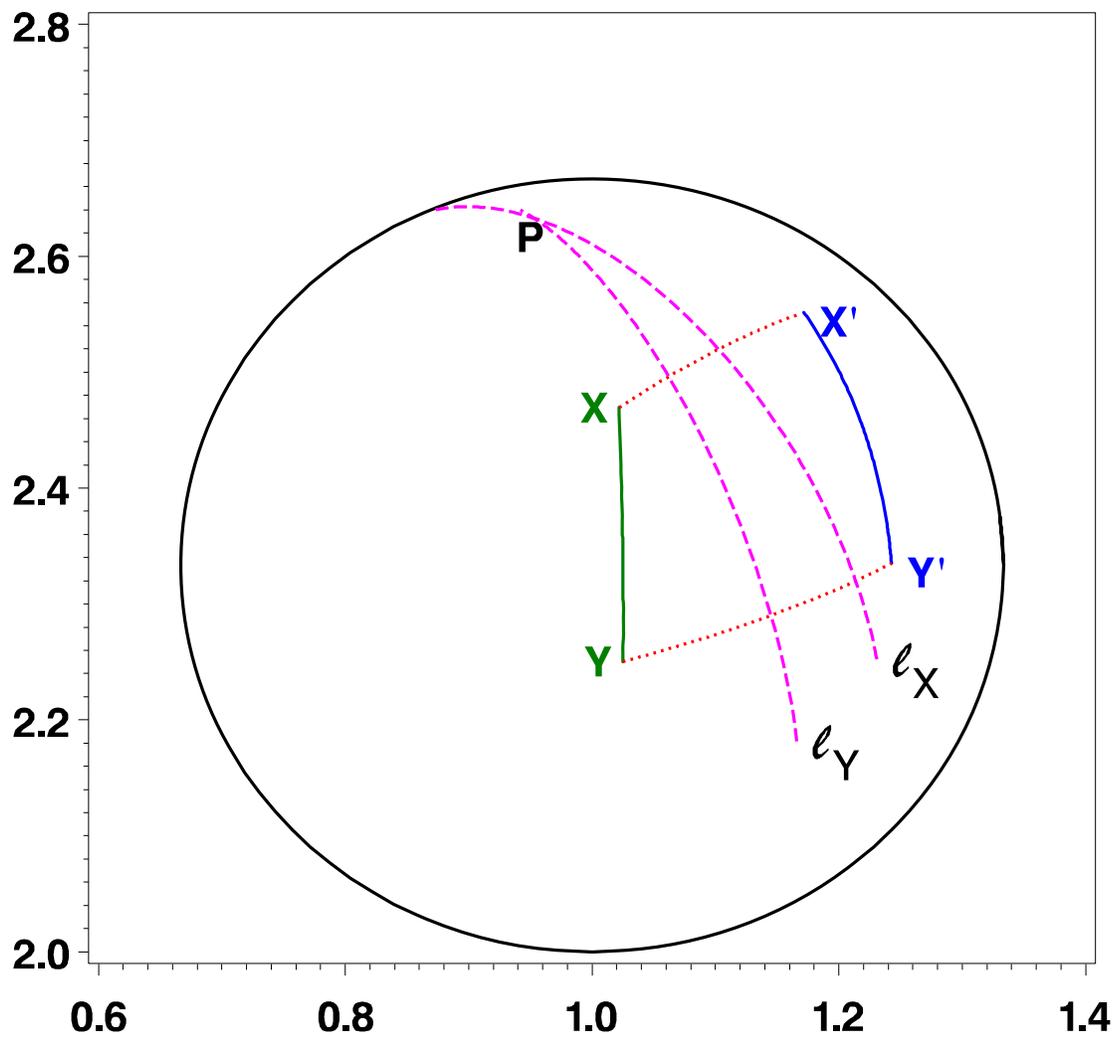

Figure 5: Geometric Location of Rotation Point P for an Isometry of the Sphere Moving XY to X'Y'